\documentclass[12pt,leqno]{article}
\begin{document}
\newtheorem{theorem}{Theorem}[section]
\newtheorem{prop}{Proposition}[section]
\title{Fedosov Quantization on Symplectic Ringed Spaces}
\author{{\small by}\vspace{2mm}\\Izu Vaisman}
\date{}
\maketitle
{\def\thefootnote{*}\footnotetext[1]%
{{\it 2000 Mathematics Subject Classification}
53D20. \newline\indent{\it Key words and phrases}:
Herz-Reinhart-Lie ringed space.
Symplectic ringed space. Fedosov quantization.}}
\begin{center} \begin{minipage}{12cm}
A{\footnotesize BSTRACT. We expose the basics of the Fedosov
quantization procedure, placed in the general framework of
symplectic ringed spaces. This framework also includes some
Poisson manifolds with non regular Poisson structures,
presymplectic manifolds, complex analytic symplectic manifolds,
etc.} \end{minipage} \end{center} \vspace{5mm}\noindent On a
symplectic manifold $(M,\omega)$ {\it Fedosov quantization} is an
embedding of the algebra $C^\infty(M,{\mathbf C})[[h]]$ of formal
power series in $h$, with complex valued differentiable functions
on $M$ as coefficients, into the algebra of the cross sections of
the {\it Weyl algebras bundle} $W(TM)$ by means of parallel
translation with respect to a {\it generalized Abelian
connection}. This paper is an exposition of the basics of Fedosov
quantization. The difference between our exposition and that of
the original works \cite{F1,F2} consists in the fact that we place
Fedosov's construction in the general framework of {\it symplectic
ringed spaces}. The generalization is purely formal, and should be
seen as folklore, but, it allows for new applications including a
class of Poisson manifolds with possibly non regular Poisson
structure, presymplectic manifolds, holomorphic symplectic
manifolds, etc. On the other hand, it is important to notice that
in the general case one has an obstruction to the existence of a
connection. Therefore, Fedosov quantization can be used only on
symplectic ringed spaces where this obstruction vanishes.
\section{Symplectic Ringed Spaces}
A {\it Herz-Reinhart Lie algebra}, shortly {\it HRL-algebra}, (or
{\it Lie pseudo-algebra} \cite{{H},{Mz}}) $L$ over a pair $(K,C)$,
where $K$ is a commutative ring with unit and $C$ is a commutative
$K$-algebra with unit, is a Lie algebra over $K$ which also is a
$C$-module, and is endowed with a mapping $\iota:L\rightarrow
{\mathcal D}$, the algebra of derivations of $C$ over $K$, that is
both a $K$-Lie algebra and a $C$-module homomorphism, with the
compatibility condition
$$[X,fY]=f[X,Y]+(\iota(X)f)Y\hspace{5mm}(f\in C;\:X,Y\in L).\leqno{(1.1)}$$

The fundamental example is $K=\mathbf{R}$, $C=C^\infty(M)$,
$L=\Gamma TM$, where $M$ is an arbitrary differentiable manifold,
and $\Gamma$ always denotes spaces of cross sections of vector
bundles or sheaves.

Now, we define a {\it Herz-Reinhart-Lie (HRL)-ringed space}, as a
topological space $M$ endowed with a sheaf ${\mathcal C}$ of
commutative $K$-algebras with unit and a {\it sheaf ${\mathcal L}$
of HRL-algebras} over $(K,{\mathcal C})$. This latter notion has
the obvious definition namely, for each open subset $U$ of $M$ the
space of sections $\Gamma_U(\mathcal L)$ is an HRL-algebra over
$(K,\Gamma_U(\mathcal C))$, and the restrictions are homomorphisms
of HRL-algebras. The sheaves of HRL-algebras were studied in
\cite{KT} under the name of sheaves of {\it twisted Lie algebras}.
An HRL-ringed space whose sheaf ${\mathcal L}$ is a locally free
$\mathcal C$-module of rank $m$ will be called an HRL-{\it ringed
space of rank} $m$.
\vspace{2mm}\\
{\bf Example 1.1} An $m$-dimensional differentiable (respectively,
complex analytic) manifold $M$ is an HRL-ringed space of rank $m$,
with $K={\mathbf R}$ (respectively, $K={\mathbf C}$), ${\mathcal
C}={\mathcal C}^\infty(M)$ (respectively $\mathcal C^{{\rm
hol}}(M)$), the sheaf of germs of differentiable (respectively,
holomorphic) functions on $M$, and ${\mathcal L}$ the sheaf of
germs of differentiable (respectively, holomorphic) vector fields
on $M$.
\vspace{2mm}\\
{\bf Example 1.2} Let $p:A \rightarrow M$ be a differentiable
vector bundle of rank $m$ which is a Lie algebroid of anchor
$\iota:A\rightarrow TM$ \cite{Mz}. Then, if we take $K={\mathbf
R}$, ${\mathcal C}=\mathcal C^\infty(M)$, ${\mathcal L}$ the sheaf
$\underline{A}$ of germs of differentiable cross sections of $A$,
and $\iota$ induced by the anchor, we again get an HRL-ringed
space of rank $m$.
\vspace{2mm}\\
{\bf Example 1.3} Let $M$ be a differentiable manifold endowed
with a foliation ${\mathcal F}$ of codimension $m$, and let
$\nu{\mathcal F}$ be the transversal bundle of ${\mathcal F}$.
Take $K={\mathbf R}$, ${\mathcal C}$ the sheaf of germs of
differentiable functions on $M$ which are constant along the
leaves of ${\mathcal F}$ ({\it foliated functions}), ${\mathcal
L}$ the sheaf of germs of foliated cross sections of $\nu{\mathcal
F}$, and $\iota(X)f=\bar Xf$, where $f\in{\mathcal C},
X\in{\mathcal L}$, and $\bar X$ is a germ of foliated vector field
on $M$ which projects onto $X$. The result is a structure of
HRL-ringed space of rank $m$. (We refer to \cite{M} for the theory
of foliated manifolds.)\vspace{2mm}
\par
Because the structure of HRL-ringed space is similar to that of
Lie algebroid, the same formulas as in the latter case
\cite{{KS1},{Mz}} yield a differential calculus for HRL-ringed
spaces.

If we refer to the sheaf $$\Omega^k(M) :=Alt_{\mathcal
C}({\mathcal L}^k,{\mathcal C})\leqno{(1.2)}$$ as the sheaf of
germs of differential $k$-forms ($:=$ denotes a definition), there
exists an exterior differential
$d:\Omega^k(M)\rightarrow\Omega^{k+1}(M)$ defined for
$\lambda\in\Omega^k(M)$ by
$$d\lambda(X_0,\ldots,X_k)=\sum_{i=0}^k
(-1)^i\iota(X_i)(\lambda(X_0,\ldots,\hat
X_i,\ldots,X_k))\leqno{(1.3)}$$
$$+\sum_{i<j=1}^{k}(-1)^{i+j}\lambda([X_i,X_j],X_0,\ldots,
\hat X_i,\ldots,\hat X_j,\ldots,X_k).$$ In (1.3), and in all the
similar formulas of this paper, the germs always are at the same
point of $M$. The operator $d$ satisfies $d^2=0$, and leads to a
usual definition of de Rham cohomology spaces $H^k(\mathcal L)$.
Notice also the existence of the {\it wedge product} which makes
the graded exterior algebra $\{\Omega^k(M)\}_{k\in{\mathbf N}}$
into a differential graded algebra with respect to the exterior
differential $d$.

Furthermore, if the sheaf
$\mathcal M$ over $M$ is a $\mathcal C$-module, we define the spaces of
${\mathcal M}$-valued differential forms by
$$\Omega^k(M,{\mathcal M})
:=Alt_{\mathcal C}{\mathcal L}^k,{\mathcal M}),\leqno{(1.4)}$$
and a {\it connection} on $\mathcal M$ is a $K$-linear homomorphism
$$\nabla:{\mathcal M}\rightarrow\Omega^1(M,{\mathcal M}),\leqno{(1.5)}$$
which satisfies the condition
$$\nabla(f\xi)=f\nabla \xi+(df)\xi\hspace{5mm}(f\in{\mathcal
C},\,\xi\in{\mathcal M}).\leqno{(1.6)}$$ We will also use the
notation  $\nabla_Y \xi:=(\nabla \xi)(Y)$.

Using a connection, and in analogy with formula (1.3), for all
$\lambda\in\Omega^k(M,{\mathcal M})$ one defines
$$\nabla\lambda(X_0,\ldots,X_k)=
\sum_{i=0}^k(-1)^i\nabla_{X_i}(\lambda(X_0,\ldots,\hat
X_i,\ldots,X_k))\leqno{(1.7)}$$
$$+\sum_{i<j=1}^{k}(-1)^{i+j}\lambda([X_i,X_j],X_0,\ldots,
\hat X_i,\ldots,\hat X_j,\ldots,X_k)\in\Omega^{k+1}
(M,\mathcal M).$$

In particular, (1.7) may be used to define the {\it curvature of
the connection} $\nabla$
$$\Phi=\nabla^2:{\mathcal M}\rightarrow\Omega^2(M,{\mathcal
M}),\leqno{(1.8)}$$ where
$$(\Phi\xi)(Y,Z)=\nabla_Y\nabla_Z \xi-\nabla_Z\nabla_Y
\xi-\nabla_{[Y,Z]}\xi. \leqno{(1.9)}$$
It is easy to check that
$\Phi$ is a homomorphism over ${\mathcal C}$, and satisfies the
{\it Bianchi identity}
$$\nabla(\Phi(\xi))(Y,Z,U)=\sum_{Cycl.(Y,Z,U)}\Phi(\nabla_Y{\xi})(Z,U).
\leqno{(1.10)}$$ In (1.9) and (1.10) $Y,Z,U\in{\mathcal L}$, and
$\xi\in \mathcal M$.

The operators $\nabla_Y$ can be extended to any {\it tensorial
sheaf} produced from $\mathcal M$ by the usual formulas of
differential geometry. In the whole paper, tensor and wedge
product sheaves are seen as sheaves of $\mathcal C$-multilinear
morphisms (e.g., \cite{V1}). In particular, the curvature $\Phi$
can also be seen as a $2$-{\it form with values in} $Hom_{\mathcal
C}({\mathcal M},{\mathcal M})$, and, then, the Bianchi identity
takes the simple classical form \cite{{H},{V1}}
$$\nabla\Phi=0.\leqno{(1.11)}$$

Our interest will be in connections on $\mathcal L$, also called
{\it connections on the HRL-ringed space} $M$. In this case
another important invariant is the {\it torsion}, defined
by $$T(X,Y)=\nabla_X Y-\nabla_Y X-[X,Y],\hspace{5mm}(X,Y\in\mathcal
L). \leqno{(1.12)}$$
From (1.1) and (1.6), it follows that $T\in\Omega^2(M,\mathcal L)$, and a
simple computation yields the {\it torsion Bianchi identity}
$$(\nabla T)(X,Y,Z)=\sum_{Cycl.(X,Y,Z)}\Phi(X)(Y,Z))\hspace{5mm}
(X,Y,Z\in\Gamma\mathcal L).\leqno{(1.13)}$$

The torsion of a connection on $\mathcal L$ yields the following
expression of the exterior differential (1.3)
$$d\lambda(X_0,\ldots,X_k)=\sum_{i=0}^k
(-1)^i(\nabla_{X_i}\lambda)(X_0,\ldots,\hat
X_i,\ldots,X_k)\leqno{(1.14)}$$
$$+\sum_{i<j=1}^{k}(-1)^{i+j}\lambda(T(X_i,X_j),X_0,\ldots,
\hat X_i,\ldots,\hat X_j,\ldots,X_k).$$
This formula suggests considering the operator
$$\nabla\lambda(X_0,\ldots,X_k)=\sum_{i=0}^k
(-1)^i(\nabla_{{X_i}}\lambda)(X_0,\ldots,\hat
X_i,\ldots,X_k)),\leqno{(1.15)}$$
which coincides with $d\lambda$ if the torsion of $\nabla$ is zero.

From (1.6), we see that the difference $D=\nabla^2-\nabla^1$ of
two connections on $M$ is a global section of $Hom_{\mathcal
C}(\mathcal L,\Omega^1 (M,\mathcal L)=Hom_{\mathcal C}(\mathcal
L\times\mathcal L,\mathcal L)$, i.e., a {\it tensor field}. Using
this remark it follows easily that if $\nabla$ is an arbitrary
connection, then
$$\nabla^0_X Y:=\frac{1}{2}(\nabla_X Y+\nabla_Y X+[X,Y]) \leqno{(1.16)}$$
is a torsionless connection.

If $\mathcal M$ is a locally free $\mathcal C$-module of finite
rank $s$ over a ringed space $(M,\mathcal C,\mathcal L)$ of rank
$m$, we may use local bases $(X_{i})_{i=1}^{m}$ of $\mathcal L$,
and local bases $(\xi_{u})_{u=1}^{s}$ of $\mathcal M$, and get the
local equations of the connection $\nabla$
$$\nabla_{X_{i}}\xi_u=\sum_{v=1}^{s}\Gamma_{iu}^{v}\xi_{v},\leqno{(1.17)}$$
where the {\it connection coefficients} $\Gamma_{iu}^{v}$ are local
sections of $\mathcal C$. The local equations (1.17) may be used in
exactly the same way as in classical differential geometry e.g.,
\cite{V1}. For instance, let us take the case $\mathcal M=\mathcal L$, and
look for the local expression of the torsion. For this purpose, we first
notice the existence of {\it structure equations}:
$$[X_{i},X_{j}]=\sum_{k=1}^{m}s^{k}_{ij}X_{k},\hspace{5mm}(s^{k}_{ij}=-
s^{k}_{ji}\in\mathcal C).$$ Then, if we put
$T(X_{i},X_{j})=\sum_{k=1}^{m}T^{k}_{ij}X_{k}$, we get
$$T^{k}_{ij}=\Gamma^{k}_{ij}-\Gamma^{k}_{ji}-s^{k}_{ij}.$$

Clearly, in the case of a differentiable manifold $M$ a connection
(1.5) on $\mathcal M=\mathcal L$ is just a linear connection on
$M$. In the case of a complex analytic manifold a connection is
the same thing as a holomorphic connection on the manifold, and it
exists iff the {\it Atiyah class} of the manifold vanishes
\cite{A}. In the case of the transversal bundle of a foliation a
connection (1.5) on $\mathcal L$ is defined by a {\it projectable
connection} of the foliation, and it exists iff the {\it Atiyah
class of the foliation} vanishes \cite{M}. Therefore, not every
ringed space has connections. In the case of a Lie algebroid a
connection (1.5) is a {\it connection of the algebroid}, and any
Lie algebroid has connections.

As a matter of fact, the Atiyah class method \cite{H,KT}, yields a general result
namely,
\begin{theorem} For each HRL-ringed space of finite rank $(M,K,\mathcal
C,\mathcal L)$, there exists a well defined sheaf-cohomology class
$a(M)\in H^1(M,Hom_{\mathcal C}(\mathcal L, Hom_{\mathcal
C}(\mathcal L$, $\mathcal L))$ such that existence of a connection on
$\mathcal L$ is equivalent with $a(M)=0$.\end{theorem} \noindent{\bf
Proof.} For any HRL-ringed space $(M,K,\mathcal C,\mathcal L)$
there exists a sheaf $\mathcal J$, which plays the role of the jet
bundle $J^1 TM$ of a differentiable manifold $M$. Namely,
$\mathcal J$ is the submodule of $\mathcal L\oplus
Hom_{K}(\mathcal L,\mathcal L)$ which consists of the pairs
$(X,-{\rm ad}_{X}+\varphi)$, where $X\in\mathcal L,\varphi\in
Hom_{\mathcal C}(\mathcal L,\mathcal L)$ and ${\rm
ad}_X:=[X,\;\;]$. The fact that $\mathcal J$ is a $\mathcal
C$-module follows from
$$(fX,-f{\rm ad}_{X}+f\varphi)=(fX,-{\rm ad}_{fX}+f\varphi-(df)X).
\leqno{(1.18)}$$
Now, we notice the existence of the following exact sequence of
${\mathcal C}$-module sheaves
$$0\rightarrow\mathcal F:=Hom_{\mathcal C}(\mathcal L,\mathcal L)
\stackrel{\epsilon}{\rightarrow}\mathcal J
\stackrel{p}\rightarrow\mathcal L\rightarrow0,\leqno{(1.19)}$$
where  $\epsilon(\varphi)=(0,\varphi)$,
$p(X,-{\rm ad}_X+\varphi)=X$.

Then, there exists a connection $\nabla$ on $\mathcal L$ iff there
exists a splitting of (1.19) i.e., a homomorphism $\psi:\mathcal
J\rightarrow{\mathcal F}$ such that $\psi\circ\epsilon=id$.
Indeed, if $\psi$ is given, $\nabla X= \psi(X,-{\rm ad}_X)$ is a
connection. Conversely, if $\nabla$ is a connection, $\psi(X,-{\rm
ad}_X+\varphi)= \nabla X+\varphi$ is the required homomorphism.

Furthermore, if $\mathcal L$ is locally free and of finite rank,
the sequence (1.19) behaves as a sequence of finite dimensional
vector spaces, and leads to the exact sequence
$$0\rightarrow Hom(\mathcal L,\mathcal F)\stackrel{p'}{\rightarrow}
Hom(\mathcal J,\mathcal F)\stackrel{\epsilon'}{\rightarrow}
Hom(\mathcal F,\mathcal F)\rightarrow0,\leqno{(1.20)}$$ then to
the corresponding exact sequence of sheaf-cohomology
$$0\rightarrow H^0(M,\mathcal A)\stackrel{p^*}{\rightarrow}
H^0(M,\mathcal B)\stackrel{\epsilon^*}{\rightarrow} H^0(M,\mathcal
E)\stackrel{\delta}{\rightarrow} H^1(M,\mathcal
A)\rightarrow\cdots \:.\leqno{(1.21)}$$ In (1.20) and (1.21), the
morphisms $\epsilon', \epsilon^*,p',p^*$ are induced by $\epsilon$
and $p$ of (1.19), $\mathcal A,\mathcal B,\mathcal E$ are the
second, third and fourth sheaf of the sequence (1.20),
respectively, and $\delta$ is the connecting morphism.

It follows that the condition for the existence of a connection is
that the identity belongs to the image of $\epsilon^*$, and this
is equivalent to $\delta(id)=0$. Hence, if we define the {\it
Atiyah class} by $a(M)=\delta(id) \in H^1(M,Hom(\mathcal
L,\mathcal F))$, we precisely have the required result. Q.e.d.
\vspace{2mm}\\
\indent Now, again, let $(M,\mathcal C,\mathcal L)$ be an
HRL-ringed space. A $2$-form $\omega\in\Omega^2(M)$ will be called
{\it non degenerate} if the sheaf homomorphism
$\flat_\omega:\mathcal L\rightarrow\Omega^1(M)$ defined by
$$\flat_\omega(X)(Y)=\omega(X,Y)\hspace{5mm} (X,Y\in\mathcal L)$$
is an isomorphism. The inverse of this isomorphism will be denoted
by $\sharp_{\omega}$. A $2$-form $\omega$ which is non degenerate
and {\it closed}, i.e., $d\omega=0$ is called a {\it symplectic
form}, and, then, $(M,\mathcal C,\mathcal L,\omega)$ is called a
{\it symplectic ringed space}.
\vspace{2mm}\\
{\bf Example 1.4} The differentiable and holomorphic symplectic
manifolds obviously are examples of symplectic ringed spaces.
\vspace{2mm}\\
{\bf Example 1.5} Let $(M,\omega)$, where $M$ is a
$(2n+h)$-dimensional differentiable manifold, and $\omega$ is a
closed $2$-form of rank $2n$ over $M$, be a {\it presymplectic
manifold}. It is well known (e.g., \cite{LM}) that ${\rm
ker}\,\omega$ is tangent to an $h$-dimensional foliation $\mathcal
S$, called the {\it characteristic foliation}, and that the form
$\omega$ is projectable with respect to this foliation. By looking
at the earlier Example 1.3, we see that the sheaf $\mathcal L$ of
germs of projectable cross sections of the transversal bundle
$\nu{\mathcal S}$ defines a ringed space structure over $M$
($\mathcal C=\mathcal C^\infty(M)$),  which is endowed with the
symplectic form induced by $\omega$.
\vspace{2mm}\\
{\bf Example 1.6} As in Example 1.2, let $A\rightarrow M$ be a Lie
algebroid of rank $2m$, with the anchor map $\iota: A\rightarrow
TM$. Then, any non degenerate cross section $\omega\in
\Gamma\wedge^2 A^*$ which is closed with respect to the exterior
differential $d_A$ (see \cite{{KS1},{Mz}}) makes the ringed space
$(M,\mathcal C^\infty(M),\underline {A})$ into a symplectic ringed
space. It is also convenient to say that $(A,\omega)$ is a {\it
symplectic Lie algebroid}. Following are some concrete examples of
symplectic Lie algebroids.
\vspace{2mm}\\
{\bf Example 1.7} If $P\in\Gamma\wedge^2 TM$ is a regular Poisson
structure of $M$, and if $\mathcal S$ is the symplectic foliation
of $P$ (e.g., \cite{V2}), then $T\mathcal S$ with the leafwise
$2$-form provided by the Poisson bracket of $P$ is a symplectic
Lie algebroid.
\vspace{2mm}\\
{\bf Example 1.8} Let $(M,W)$ be an arbitrary Poisson manifold.
Then $T^*M$ is a Lie algebroid of anchor $\iota=\sharp_W$. A
symplectic structure on this algebroid is a non degenerate
bivector field $Q$ on $M$ which is a cocycle in the
Poisson-Lichnerowicz cohomology i.e., $[W,Q]=0$ (e.g.,
\cite{V2}).\vspace{2mm}\\
\indent Some of the usual symplectic notions straightforwardly
transfer to symplectic ringed spaces. For instance, we may define
$\Pi\in Alt_{\mathcal C}(\Omega^1(M)\times\Omega^1(M),\mathcal C)$
by
$$\Pi(\sigma_1,\sigma_2)=<\sigma_1,\sharp_\omega\sigma_2>,\hspace{5mm}
\sigma_1,\sigma_2\in\Omega^1(M),$$ and also define
$\sharp_\Pi:=\sharp_\omega$. Furthermore, we may define the {\it
Hamiltonian gradient} $X_f$ of any $f\in\mathcal C$ by
$X_f=-\sharp_\Pi df$, and the {\it Poisson bracket}
$$\{f,g\}:=\omega(X_f,X_g)=\iota(X_f)g=-\iota(X_g)f=\Pi(df,dg).\leqno{(1.22)}$$

The evaluation
$$0=d\omega(X_f,X_g,X_h)=\sum_{Cycl.(f,g,h)}\iota(X_f)\omega(X_g,X_h)-
\sum_{Cycl.(f,g,h)}\omega([X_f,X_g],X_h)$$
$$=\sum_{Cycl.(f,g,h)}\{f,\{g,h\}\}-\sum_{Cycl.(f,g,h)}<dh,[X_f,X_g]>
=\sum_{Cycl.(f,g,h)}\{f,\{g,h\}\}$$
$$-\sum_{Cycl.(f,g,h)}(\iota[X_f,X_g])h
=\sum_{Cycl.(f,g,h)}\{f,\{g,h\}\}$$
$$-\sum_{Cycl.(f,g,h)}(\iota(X_f)\iota(X_g)-\iota(X_g)\iota(X_f))h
=-\sum_{Cycl.(f,g,h)}\{f,\{g,h\}\}$$ shows that the bracket (1.22)
satisfies the Jacobi identity. Since the Leibniz rule also
obviously holds, we have got a structure of {\it Poisson algebras
sheaf} on $\mathcal C$. (See, for instance, \cite{V2} for the
definition of a Poisson algebra.)

However, not all the classical symplectic properties hold. For
instance, the Jacobi identity for the Poisson bracket (1.22)
implies
$$X_{\{f,g\}}-[X_f,X_g]\in{\rm ker}\,\iota,\leqno{(1.23)}$$
hence, the result is zero only at the {\it injectivity points}
of $\iota$.
Another negative example is obtained if we consider the operation
of a {\it Schouten-Nijenhuis bracket} on an HRL-ringed space (e.g.,
\cite{KS1}). Then, the general algebraic computations of \cite{GD} hold,
and we have the formula
$$[\Pi,\Pi](df,dg,dh)=2\sum_{Cycl(f,g,h)}\{\{f,g\},h\}.\leqno{(1.24)}$$
But, since $\{df\;/\;f\in\mathcal C\}$ may not span $\Omega^1(M)$,
generally, we may have $[\Pi,\Pi]\neq0$.

Let us come back to the symplectic Lie algebroid $(A,\omega)$ of
Example 1.6. Then, the symplectic objects $\flat_\omega$, $\Pi$,
$X_f$, etc. have interpretations in terms of vector bundles: an
isomorphism $\flat_\omega:A\rightarrow A^*$, a cross section
$\Pi\in\Gamma\wedge^2A$, cross sections $X_f\in\Gamma A$, etc. The
Poisson bracket (1.22) becomes a Poisson algebra structure on
$C^\infty(M)$, and there exists a corresponding Poisson bivector
field $P\in\Gamma\wedge^2TM$ such that
$$\{f,g\}=\Pi(d_Af,d_Ag)=P(df,dg),\hspace{2mm}[P,P]=0,\leqno{(1.25)}$$
where the final bracket is the usual Schouten-Nijenhuis bracket on
$M$.

For a better understanding of
the relation between $P$ and $\Pi$, let us consider the transposed
homomorphism $\iota':T^*M\rightarrow A^*$ of the anchor $\iota:A
\rightarrow TM$. Then we have
$$d_Af=\iota'(df),\hspace{2mm}\sharp_P=-\iota\circ\sharp_{\Pi}\circ\iota'.
\leqno{(1.26)}$$

Notice that the injectivity points of $\iota$ are the same as the surjectivity
points of $\iota'$. Accordingly, from (1.23) and (1.24), we see that the
relations
$$X_{\{f,g\}}=[X_f,X_g],\;\;[\Pi,\Pi]_A=0\leqno{(1.27)}$$
hold at the injectivity points of $\iota$. If the set of
injectivity points of $\iota$ is dense in $M$, the equalities
(1.27) hold everywhere on $M$.

The considerations above suggest the following definition: a
Poisson structure $P$ on a differentiable manifold $M$ will be
called a {\it quasi-symplectic Poisson structure} if it is induced
by a symplectic Lie algebroid, via formula (1.25).

Example 1.7 tells us that every regular Poisson structure is
quasi-symplectic. Other quasi-symplectic Poisson structures are
provided by Example 1.8. If $(M,W)$ is a Poisson manifold which
has a non singular Poisson $2$-cocycle $Q$ (i.e.,
$Q\in\Gamma\wedge^2TM$, rank$\,Q$=dim$\,M$, $[W,Q]=0$), the
corresponding object $\Pi$ of (1.25) is a usual $2$-form on $M$,
and (1.25) defines a quasi-symplectic Poisson structure given by
$$\{f,g\}=\Pi(X_f^W,X_g^W)\hspace{5mm}(f,g\in
C^\infty(M)),\leqno{(1.28)}$$
where the arguments are
the $W$-Hamiltonian vector fields of the functions $f,g$.

For a more concrete situation of this kind, let $(M,\Pi,\phi)$ be
a {\it symplectic-Nijenhuis manifold}, with the symplectic form
$\Pi$ and the Nijenhuis tensor $\phi$ (see, for instance,
\cite{V5}). Then, $M$ has a Poisson structure $W$ defined by
$\sharp_W=\phi\circ\sharp_Q$, where $Q$ is the bivector field
given by $\sharp_Q=\flat_{\Pi}^{-1}$. $W$ is the first structure
of the Poisson hierarchy of the symplectic-Nijenhuis manifold $M$,
and it is compatible with the Poisson structure $Q$ of $M$ i.e.,
$[W,Q]=0$ \cite{V5}. Since $Q$ is non degenerate, we are in the
situation described by the previous paragraph, and we get a
quasi-symplectic Poisson structure
$$\{f,g\}=\Pi(\sharp_Wdf,\sharp_Wdg)=<\flat_\Pi\sharp_Wdf,\sharp_Wdg>$$
$$=-<\sharp_W\flat_\Pi\sharp_Wdf,dg>=-<\phi\sharp_Wdf,dg>=-<\phi^2\sharp_Qdf,dg>.$$
Up to the sign, this is the second structure of the Poisson
hierarchy of $(M,\Pi,\phi)$ \cite{V5}. Thus, we have proven
\begin{prop} The second Poisson structure of the Poisson hierarchy
of a symplectic-Nijenhuis manifold is a quasi-symplectic Poisson
structure. \end{prop} \indent Of course, the structures of
Proposition 1.1
may have singular points.\vspace{2mm}\\
\indent One of the main ingredients of Fedosov quantization is a
symplectic connection. The well known procedure of constructing
symplectic connections on symplectic manifolds, as presented for
instance in \cite{V3}, holds without modification on a symplectic
ringed space.

Namely, if $(M,\mathcal C,\mathcal L)$ is a ringed space which
possesses a non degenerate $2$-form $\omega$ ({\it almost
symplectic ringed space}), a connection $\nabla$ on $\mathcal L$
(i.e., on $M$) is said to {\it preserve} $\omega$ if
$\nabla_{X}\omega=0$, $\forall X\in\mathcal L$. If there exists a
connection $\nabla^0$ on $M$, the formula
$$\nabla_X Y=\nabla^0_X Y+\Theta(X,Y)+A(X,Y),\leqno{(1.29)}$$
where
$$\omega(\Theta(X,Y),Z)=\frac{1}{2}(\nabla^0\omega)(Y,Z),
\;\omega(A(X,Y),Z)=B(X,Y,Z),\leqno{(1.30)}$$
$\forall X,Y,Z\in\mathcal L$, and $\forall B\in Hom_\mathcal C(\mathcal
L\times\mathcal L\times\mathcal L,\mathcal C)$ which satisfies the condition
$B(X,Y,Z)=B(X,Z,Y)$, yields all the connections of $M$ which preserve
$\omega$. $\Theta$ and $A$ are well defined by (1.30) because $\omega$ is
non degenerate.

Furthermore, formula (1.14) shows that a torsionless,
$\omega$-preserving connection may exist only if $d\omega=0$. Conversely,
if we are in this latter case, and if we assume that $\nabla^0$ has zero
torsion (e.g., this $\nabla^0$ is given by applying (1.16) to the
original $\nabla^0$), it turns out that $\nabla$ of (1.29), and with
$$B(X,Y,Z)=\frac{1}{6}[(\nabla^0_Y\omega)(X,Z)+
(\nabla^0_Z\omega)(X,Y)]\leqno{(1.31)}$$ is an
$\omega$-preserving, torsionless connection. Indeed, from (1.29),
and since $d\omega=\nabla^0\omega$ (in the sense of (1.15)), we
deduce that $\omega(T_\nabla(X,Y),Z)=0$ for the chosen value
(1.31) of $B$.

By definition, a connection which preserves the symplectic form
$\omega$ and has zero torsion is called a {\it symplectic
connection} on the symplectic ringed space $(M,\mathcal C,\mathcal
L,\omega)$. Above, we saw that if $M$ has an arbitrary connection
it also has  symplectic connections, and we wrote down the
expression of one symplectic connection $\nabla$, defined by means
of (1.31). It follows that all the symplectic connections are
given by
$$\tilde\nabla_X Y=\nabla_X Y+ C(X,Y),$$
where $Q(X,Y,Z):=\omega(C(X,Y),Z)$ is symmetric in all its arguments.

In analogy with Riemannian geometry, one defines the {\it
covariant curvature tensor}  of a symplectic connection namely,
$$S(X_1,X_2,X_3,X_4)=-\omega(X_1,\Phi(X_2)(X_3,X_4)),\leqno{(1.32)}$$
where $\Phi$ is given by (1.8). This tensor has the following
symmetry properties
$$\begin{array}{l}
S(X_1,X_2,X_3,X_4)=-S(X_1,X_2,X_4,X_3),\vspace{2mm}\\
S(X_1,X_2,X_3,X_4)+S(X_1,X_3,X_4,X_2)+S(X_1,X_4,X_2,X_3)=0,\vspace{2mm}\\
S(X_1,X_2,X_3,X_4)=S(X_2,X_1,X_3,X_4).\end{array}\leqno{(1.33)}$$
The first equality is obvious, the second is the {\it Bianchi
identity} (1.13), and the third follows by expressing the
derivations via the connection in the identity
$$(\iota(X_3)\iota(X_4)-\iota(X_4)\iota(X_{3})-\iota([X_3,X_4]))
(\omega(X_1,X_2))=0,$$
where $X_1,X_2,X_3,X_4\in\mathcal L$.

Concerning examples of symplectic connections in other cases than
symplectic manifolds, we first quote the case of a presymplectic
manifold $(M,\omega)$. It is easy to see that a symplectic
connection on $M$ seen as a symplectic ringed space, as in Example
1.5, may be identified with a connection on $M$ seen as a
differentiable manifold, which preserves $\omega$, has torsion
tangent to the characteristic foliation of $\omega$, and defines a
transversal projectable connection of the same foliation. Details
on the construction of these connections can be found in
\cite{V4}.

Another interesting case is that of a symplectic structure $Q$ on
the tangent Lie algebroid $T^*M$ of a Poisson manifold $(M,W)$
(see Example 1.8). In this case, it is natural to start with a
connection $D$ on the differentiable manifold $M$ which satisfies
the condition $D_X Q=0$ ($X\in\Gamma TM$). $D$ exists since
$(M,Q^{-1})$ is an almost symplectic manifold but, generally, $D$
has a non zero torsion $T_D$. The connection $D$ yields a
connection $\nabla$ on the Lie algebroid $T^*M$ \cite{{V2},{Xu}}
by putting
$$\nabla_\alpha\beta=D_{\sharp_W\alpha}\beta,\hspace{1cm}\alpha,
\beta\in\Gamma T^*M.\leqno{(1.34)}$$ A straightforward
computation, which uses the Lie bracket of $T^*M$
\cite{{KS1},{V2}}
$$\{\alpha,\beta\}=
L_{\sharp_W\alpha}\beta-L_{\sharp_W\beta}\alpha-d(W(\alpha,\beta)),
$$
yields the torsion $T_\nabla$ by the formula
$$<T_\nabla(\alpha,\beta),X>=\alpha(T_D(\sharp_W\beta,X))-
\beta(T_D(\sharp_W\alpha,X))-(D_XW)(\alpha,\beta).\leqno{(1.35)}$$

Therefore, generally,
$T_\nabla\neq0$, and we must apply the general algorithm to get a torsionless
symplectic connection on $(T^*M,Q)$.
The passage from $\nabla$ of (1.34) to a torsionless
connection is given by (1.16) i.e.,
$$\nabla^0_\alpha \beta=\frac{1}{2}(\nabla_{\alpha}\beta+
\nabla_{\beta}\alpha
+\{\alpha,\beta\}).\leqno{(1.36)}$$

Furthermore, from (1.34), (1.36), and $d_{T*M}Q=0$ we get
$$(\nabla^0_\alpha Q)(\beta,\gamma)=\frac{1}{2}Q(\alpha,T_\nabla(\beta,
\gamma)).\leqno{(1.37)}$$

Finally, from (1.30), (1.31), (1.37), we deduce the expression of
a symplectic connection
$$\check\nabla_\alpha\beta=\nabla^0_\alpha\beta+
E(\alpha,\beta),\leqno{(1.38)}$$
where
$$Q(E(\alpha,\beta),\gamma)=\frac{1}{4}Q(\alpha,T_\nabla(\beta,\gamma))
\leqno{(1.39)}$$ $$+
\frac{1}{12}[Q(\beta,T_\nabla(\alpha,\gamma))
+Q(\gamma,T_\nabla(\alpha,\beta))].$$
\section{Fedosov Quantization}
In this section, we describe the Fedosov quantization procedure,
following the original works \cite{F1,F2}, with minor
modifications, and emphasizing the ringed space setting. Other
expositions of this procedure can be found in \cite{Dz,Kr,Xu2},
etc.

First we recall the construction of the {\it Weyl algebra} of a
complex symplectic vector space $(E,\omega\in\wedge^2 E)$ of
dimension $2m$.

We begin with the associative, commutative algebra of formal
Laurent series in the parameter $h$
$$W(E):=\{w=\sum_{k=-\infty}^\infty\sum_{i=0}^\infty t_{ki}h^k
\;/\;t_{ki}\in\odot^i E^*,\;i+2k\geq0\},\leqno{(2.1)}$$ where
$\odot$ denotes the symmetric tensor product, and it defines the
multiplication in $W(E)$.

This algebra will be graded by asking
$${\rm deg}\,t_{ki}=i,\; {\rm deg}\,h=2,\leqno{(2.2)}$$
and we will write the element $w\in W(E)$ which appears in (2.1)
as
$$w=\sum_{p=0}^\infty\hat w_{p},\hspace{5mm}\hat
w_p=\sum_{2k+i=p}t_{ki}h^k.
\leqno{(2.3)}$$
In terms of vector spaces, this means
$$W(E)=\prod_{p=0}^\infty V_{p}(E),\hspace{5mm}
V_{p}(E)=\prod_{k=-\infty}^{[p/2]}\odot^{p-2k}E^*. \leqno{(2.4)}$$

Furthermore, consider the {\it contraction operators}
$$C_\omega^p:(\odot^i E^*)\otimes(\odot^j E^*)\rightarrow
\odot^{i+j-2p}E^*\leqno{(2.5)}$$
defined by $0$ if $p>{\rm min}(i,j)$, and by
$$C_\omega^p(\alpha,\beta)(e_1,\cdots,e_{i+j-2p})=
\sum_{\sigma\in S_{i+j-2p}}\sum_{a_1,\cdots,a_p=1}^m
\frac{1}{(i+j-2p)!}\leqno{(2.6)}$$
$$\cdot\alpha(b_{a_{1}^{*}},\cdots,b_{a_{p}^{*}},
e_{\sigma(1)},\cdots,e_{\sigma(i-p)})\beta(b_{a_{1}},\cdots,b_{a_{p}},
e_{\sigma(i-p+1)},\cdots,e_{\sigma(i+j-2p)})$$ if $p\leq{\rm
min}(i,j)$. In (2.6), $S$ is the symmetric group, $e_l\in E$
$(l=1,...,i+j-2p)$, and $(b_a,b_{a^*})$ $(a^*:=a+m)$ is an
arbitrary $\omega$-symplectic basis of $E$.

The contractions extend to series (2.1), and may be used to define
the {\it Moyal product}
$$w\circ w':=\sum_{q=0}^\infty\frac{1}{q!}\left(-\frac{ih}{2}\right)^q
C^q_\omega(w,w')\hspace{5mm}(w,w'\in W(E)),\leqno{(2.7)}$$ which
makes $W(E)$ into a graded (because of (2.2)), associative
algebra, called the {\it Weyl algebra} of $(E,\omega)$.

Computations in the Weyl algebra become easy if we consider an
arbitrary basis $(u_i)_{i=1}^{2m}$ of $E$, and represent an
element $t\in\odot^q E^*$ by
$$t=t_{i_1,\cdots,i_q}y^{j_1}\cdots
y^{j_q}=\sum_{|\alpha|=q}\tau_{\alpha_1\cdots\alpha_{2m}}
(y^1)^{\alpha_1}\cdots (y^{2m})^{\alpha_{2m}}=
\sum_{|\alpha|=q}\tau_\alpha y^\alpha.\leqno{(2.8)}$$ In (2.8),
$(y^i)$ are the coordinates of a generic vector of $E$ with
respect to the basis $(u_i)$, and for their indices we use (here
and subsequently) the Einstein summation convention. Furthermore,
the coefficients $t$ are symmetric, the second equality is
obtained by collecting the various factors $y^i$ into a power of
$y^i$, and the third equality is the formal notation of its left
hand side, i.e.,
$$\alpha=(\alpha_1,\cdots,\alpha_{2m}),\;(\alpha_i\geq0),\;
y^\alpha=(y^1)^{\alpha_1}\cdots (y^{2r})^{\alpha_{2m}},\;|\alpha|=
\alpha_1+\cdots+\alpha_{2m}.$$

With this representation, an element $w\in W(E)$ becomes
$$w(y,h)=\sum_{k=-\infty}^{k=\infty}\sum_
{|\alpha|=0}^{\infty}h^k
t_{k,\alpha}y^\alpha,\hspace{5mm}|\alpha|+2k\geq0, \leqno{(2.9)}$$
the usual product of polynomials corresponds to the symmetric
tensor product, and the Moyal product is
$$w\circ w'=\sum_{q=0}^\infty\frac{1}{q!}\left(-\frac{ih}{2}\right)^q
\omega^{i_1j_1}\cdots\omega^{i_qj_q}
\frac{\partial^q w}{\partial y^{i_1}\ldots\partial y^{i_q}}
\frac{\partial^q w'}{\partial y^{j_1}\ldots\partial y^{j_q}}\leqno{(2.10)}$$
$$=exp\left(-\frac{ih}{2}\omega^{ij}\frac{\partial}{\partial z^i}
\frac{\partial}{\partial u^j}\right)(w(z,h)w'(s,h)/_{z=s=y},$$
where $\omega^{ih}\omega_{hk}=\delta_k^i$ and $\omega_{hk}$ are
the $u$-components of $\omega$.

It follows easily that the center $Z(W(E))$ is the
algebra of formal power series ${\mathbf C}[[h]]$ \cite{{F1},{F2}}.

The next step consists of {\it enlarging} the Weyl algebra $W(E)$
to the associative algebra
$$\hat W(E):=W(E)\otimes(\oplus_{q=0}^{2m}\wedge^q E^*).\leqno{(2.11)}$$
If $(u_i)$ is the basis used in (2.8), and if $(\nu^i)$ is its dual
cobasis, $\lambda\in \hat W(E)$ may be seen as
$$\lambda=
\sum_{k=-\infty}^\infty\sum_{p=0}^{\infty} \sum_{q=0}^{2m}
h^k\lambda_{kpqi_1\cdots i_pj_1\cdots j_q} y^{i_1}\cdots
y^{i_p}\nu^{j_1}\wedge\cdots\wedge\nu^{j_q},
\hspace{3mm}p+2k\geq0,\leqno{(2.12)}$$ where the coefficients are
symmetric in the indices $i$, and skew-symmetric in the indices
$j$. The product of $\hat W(E)$ is defined by (2.10) with a wedge
product of the partial derivatives which appear in that formula.

In $\hat W(E)$, the {\it commutant} is defined as the natural
extension of
$$[\lambda,\mu]=\lambda\circ\mu-(-1)^{\tilde\lambda\tilde\mu}\mu\circ
\lambda,\leqno{(2.13)}$$
where tilde denotes the degree of the wedge product factor of an element
of $\hat W(E)$.

The center $Z(\hat W(E))$ is ${\mathbf
C}[[h]]\otimes(\oplus_{q=0}^{2m} \wedge^q E^*)$, and one has the
{\it central projections} \cite{F1,F2} of $\lambda$ of (2.12):
$$\lambda_0:=\sum_{k=0}^{\infty}\sum_{q=0}^{2m}h^k
\lambda_{k0qj_1\cdots j_q}\nu^{j_1}\wedge\cdots\wedge
\nu^{j_q}\in Z(\hat W(E)),\leqno{(2.14)}$$
$$\lambda_{00}:=\sum_{k=0}^{\infty}h^k\lambda_{k00}\in Z(W(E)).\leqno{(2.15)}$$

The following basis-independent operators are essential in the
subsequent computations \cite{{F1},{F2}}:
$$\delta\lambda:=\sum_{j=1}^{2m}\nu^j\wedge\frac{\partial\lambda}{\partial y^j},
\;\delta^*\lambda=\sum_{j=1}^{2m}y^j
(i(u_{j})\lambda).\leqno{(2.16)}$$ These operators satisfy the
properties
$$\delta^2=0,\;\delta^{*2}=0,\;\delta(a\circ b)=(\delta a)\circ b+
(-1)^{\tilde a}a\circ\delta b.\leqno{(2.17)}$$
Furthermore, $\forall\lambda\in\hat W(E)$,
one can check the {\it Hodge decomposition} formula
$$\lambda=\delta\delta^{-1}\lambda+\delta^{-1}\delta\lambda+\lambda_{00},
\leqno{(2.18)}$$
where $\delta^{-1}$ is defined on the $(p,q)$-term of (2.12) by
$$\delta^{-1}\lambda:=\frac{1}{p+q}\delta^*\lambda.\leqno{(2.19)}$$

We intend to apply the previous algebraic constructions to
symplectic ringed spaces. To be able to do so, in what follows we
assume that $\mathcal C$ is a subsheaf of the sheaf of germs of
continuous, complex valued functions on $M$ and that $(M,\mathcal
C,\mathcal L,\omega)$ is a symplectic ringed space of finite rank.
Then $\mathcal L$ must the sheaf of germs of cross sections of a
complex symplectic vector bundle over $M$ \cite{Ten}. This implies
that the rank of the space is even, say $2m$, and that $\mathcal
L$ has local {\it symplectic bases} i.e., local bases
$(X_i)_{i=1}^{2m}$ such that
$$\omega(X_a,X_b)=0,\;\omega(X_a,X_{b+m})=\delta_{ab},\;
\omega(X_{a+m},X_{b+m})=0,\hspace{3mm}a,b=1,\ldots r.$$

Then the previous algebraic constructions may be performed on each
stalk of $\mathcal L$, and using local symplectic bases of germs
$(X_i)$. Accordingly, we get {\it sheaves of Weyl algebras}
$\mathcal W(\mathcal L)$, $\hat\mathcal W(\mathcal L)$, and the
formulas developed earlier in this section hold, with germs
instead of algebraic tensors overall. In particular, there is a
{\it central sheaf} $\mathcal Z(\mathcal W(\mathcal L))=\mathcal
C[[h]]$, which consists of germs of formal power series in $h$,
and a {\it central sheaf} $\mathcal Z(\hat\mathcal W(\mathcal
L))=\mathcal C[[h]]\otimes (\oplus_{q=0}^{2m}\wedge^q\mathcal L)$.

In the ringed setting, Fedosov's quantization  will be an
embedding of $\mathcal C[[h]]$ onto the {\it parallel germs} of
$\mathcal W(\mathcal L)$ with respect to a {\it generalized
symplectic connection}. Accordingly, we shall assume that
$(M,\mathcal C,\mathcal L,\omega)$ has connections (it has a
vanishing Atiyah class $a(M)$, as defined in Section 1), and take
a torsionless, symplectic connection $\nabla$ on this space. Then,
$\nabla$ extends to a {\it covariant exterior differential}
$\nabla:\hat\mathcal W(\mathcal L)\rightarrow\hat\mathcal
W(\mathcal L)$ defined by \cite{{F1},{F2}}
$$\nabla(t\otimes\theta)=\sum_{i=1}^{2m}\nu^i\wedge\nabla_{X_i}
(t\otimes\theta)\hspace{5mm}(t\in\odot^q\mathcal
L^*,\theta\in\wedge^s \mathcal L^*),\leqno{(2.20)}$$ where $(X_i)$
is a local basis of $\mathcal L$ and $(\nu^i)$ is the dual
cobasis. The definition is invariant by a change of the local
basis, and (2.20) reduces to (1.15) if there is no symmetric
factor $t$. Since the connection has zero torsion, in the case of
a differential form $\theta$ one has $\nabla\theta=d\theta$.
Furthermore, since $\nabla\omega=0$, (2.10) shows that
$$\nabla(\lambda\circ\mu)=(\nabla\lambda)\circ\mu+(-1)^{\tilde\lambda}
\lambda\circ\nabla\mu, \hspace{5mm}
\lambda,\mu\in\hat\mathcal W(\mathcal L).\leqno{(2.21)}$$

The operator (2.20) is what is actually needed in Fedosov
quantization, and, subsequently, we will think of this operator
when referring to a connection. Fedosov \cite{{F1},{F2}} writes
the operator $\nabla$ in a convenient way as follows. Consider the
local equations (1.17) of the connection $\nabla$, and assume that
the basis $(\xi_i=X_i)$ used in these equations is symplectic.
Then, the symplectic character of the connection is equivalent to
$$\Gamma_{ijk}=\Gamma_{jik}\hspace{5mm}(\Gamma_{ijk}:=\omega_{is}\Gamma_{jk}^s).
\leqno{(2.22)}$$
Accordingly, $\forall x\in M$, there exists a germ
$$\Gamma:=\frac{1}{2}\Gamma_{ijk}y^i y^j\nu^k\in(\hat\mathcal W(\mathcal L))_x,
\leqno{(2.23)}$$ and {\it Fedosov's formula} is
$$\nabla\lambda=d\lambda+\frac{i}{h}[\Gamma,\lambda],\leqno{(2.24)}$$
for $\lambda$ given by (2.12), and with $d$ applied as if $h,y$ would be
constants. Formula (2.24) is easily checked for $\lambda=\lambda_i \nu^i$
and $\lambda=\lambda_i y^i$, and it holds in the general case because
$d$ and the commutant $[\Gamma,\;\;]$ are derivations of
$\hat\mathcal W(\mathcal L)$. Notice that the germs (2.23) do not define
a global section of $\hat\mathcal W(\mathcal L)$.

The same method yields the formulas \cite{{F1},{F2}}
$$\delta\lambda=-\frac{i}{h}
[\varpi,\lambda],\hspace{2mm}
\varpi:=\delta^*\omega=\omega_{ij}y^i\nu^j,\leqno{(2.25)}$$
$$\nabla\delta+\delta\nabla=0,\leqno{(2.26)}$$
$$\nabla^2\lambda=\frac{i}{h}[S,\lambda],\hspace{5mm}
S:=-\frac{1}{4}S_{ijkl}y^i y^j \nu^k\wedge\nu^l,\leqno{(2.27)}$$
$S_{ijkl}$ being the components of the covariant curvature tensor
of the symplectic connection $\nabla$, which is known to be
symmetric in the first two arguments and skew symmetric in the
last two arguments (see (1.33)).

Fedosov's formula (2.24) suggests a definition of {\it generalized symplectic
connections} \cite{{F1},{F2}} as operators
$$\hat\nabla\lambda=\nabla\lambda+\frac{i}{h}[\gamma,\lambda],\leqno{(2.28)}$$
where $$\gamma=\sum_{k=-\infty}^\infty
\sum_{p=0}^\infty h^k\gamma_{ki_1\cdots
i_pj}
y^{i_1}\cdots y^{i_p}\nu^j\hspace{3mm}(p+2k\geq0)
\leqno{(2.29)}$$ are germs which define a global section
of $\Gamma(\mathcal W(\mathcal L)\otimes\mathcal L^*)$.
Then, a straightforward computation yields the {\it generalized curvature}
$\Phi$ defined by means of the formulas
$$\hat\nabla^2\lambda=\frac{i}{h}[\Phi,\lambda],\hspace{5mm}
\Phi=S+\nabla\gamma+\frac{i}{h}\gamma^2,\leqno{(2.30)}$$
and the {\it Bianchi identity}
$$\hat\nabla\Phi=\nabla\Phi+\frac{i}{h}[\gamma,\Phi]=0.\leqno{(2.31)}$$

By definition, if $\hat\nabla^2\lambda=0$ for all $\lambda\in
\hat\mathcal W(\mathcal L)$, $\hat\nabla$ is called an {\it Abelian connection}.

Now, we come to the result which is at the heart of Fedosov quantization
\cite{{F1},{F2},{Dz},{Xu2}}\\
\begin{theorem} Let $(M,\mathcal C,\mathcal L,\omega)$ be a
symplectic ringed space of finite rank $2m$, with the sheaf
$\mathcal C$ being a subsheaf of germs of continuous, complex
valued functions, and which has connections. Then, there exist
generalized, symplectic, Abelian connections $\hat\nabla$ on $M$.
Furthermore, for any Abelian connection $\hat\nabla$, for any
$a\in\mathcal C[[h]]$, there exists a unique $\lambda\in\mathcal
W(\mathcal L)$ with central projection $\lambda_{00}=a$, such that
$\hat\nabla\lambda=0$.
\end{theorem}
\noindent{\bf Proof.}
Consider a generalized connection (2.28), where $\gamma=\varpi+r$ for
$\varpi$ given by (2.25), and for an element
$$r=\sum_{p\geq2}^\infty \hat r_{p}\in\mathcal W(\mathcal L)
\otimes\mathcal L^*\leqno{(2.32)}$$
which satisfies the condition $\delta^{-1}r=0$. The last condition
implies that $r$ has the  central projection $r_0=0$, therefore, by the
Hodge decomposition (2.18) $r=\delta^{-1}\delta r$.

The curvature form of this connection $\hat\nabla$ is given by
(2.30) and, in view of (2.10), (2.13), (2.25) and
$\nabla\omega=0$, it becomes
$$\Phi=S+\nabla r-\delta r+\frac{i}{h}r^2-\omega,\leqno{(2.33)}$$
Since $\omega\in\mathcal Z(\hat\mathcal W(\mathcal L))$, the
condition $$\delta r=S+\nabla r+\frac{i}{h}r^2\leqno{(2.34)}$$
ensures $\hat\nabla^2=0$.

We show that (2.34) has a unique solution $r$ with the required
properties. Uniqueness will ensure that the germs $r$ define a
global cross section of $\mathcal W(\mathcal L)\otimes \mathcal
L^*$. Indeed, by applying to (2.34) the operator $\delta^{-1}$, we
get
$$\delta^{-1}\delta r=r =
\delta^{-1}S+\delta^{-1}\nabla r+\frac{i}{h}\delta^{-1}
(r^2),\leqno{(2.35)}$$ which is equivalent to the recurrence
formula
$$\hat r_p=(\widehat{\delta^{-1} S})_p+
\delta^{-1}\nabla\hat r_{p-1}+\frac{i}{h}\delta^{-1}
(\sum_{i=2}^{p-3}\hat r_i\circ\hat r_{p-1-i}).\leqno{(2.36)}$$

(2.32) and (2.36) imply $\hat r_2=0$, and, then,
$$\hat r_p=(\delta^{-1}\nabla)^{p-3}\delta^{-1}S\leqno{(2.37)}$$
$$+\frac{i}{h}\sum_{s=0}^{p-7}[(\delta^{-1}\nabla)^s\delta^{-1}
(\sum_{j=3}^{p-s-2}\hat r_j\circ\hat r_{p-s-j-1})]\;(p\geq3).$$

Now, it remains to establish that the obtained germ $r$ satisfies
equation (2.34). From (2.35),
using the Hodge decomposition (2.18) and the Bianchi
identity (1.33), we get
$$\delta r=\delta\delta^{-1}S+\delta\delta^{-1}\nabla r+\frac{i}{h}
\delta\delta^{-1}r^2
=S+\nabla r+\frac{i}{h}r^2-\delta^{-1}\delta(\nabla r
+\frac{i}{h}r^2).\leqno{(2.38)}$$
Formulas (2.33), (2.38) give us the curvature
$$\Phi=\delta^{-1}\delta(\nabla r+\frac{i}{h}r^2)-\omega,$$
and, since
$$\hat\nabla \omega\stackrel{(2.28)}{=}\nabla\omega+
\frac{i}{h}[\varpi+r,\omega]=\nabla\omega=d\omega=0,$$
the Bianchi identity (2.31) yields
$$\hat\nabla\delta^{-1}\delta(\nabla r+\frac{i}{h}r^2)=0.\leqno{(2.39)}$$

But, if we look at any $\lambda$ such that $\delta\lambda=0$ and
$\lambda_{00}=0$, we have $\lambda=\delta\delta^{-1}\lambda$, and
we see that
$$\hat\nabla\delta^{-1}\lambda\stackrel{(2.28)}{=}
\nabla\delta^{-1}\lambda-\delta\delta^{-1}\lambda+\frac{i}{h}
[r,\delta^{-1}\lambda]=0$$ implies
$$\lambda=\nabla\delta^{-1}\lambda+\frac{i}{h}[r,\delta^{-1}\lambda].
\leqno{(2.40)}$$ This is equivalent to a recurrence relation
$$\hat \lambda_p=\nabla\delta^{-1}\hat \lambda_{p-1}+\frac{i}{h}
({\rm terms\;in\;}\hat \lambda_{i},\;i\leq p-4),\leqno{(2.41)}$$
which yields $\lambda=0$. Since $$\lambda=\delta(\nabla
r+\frac{i}{h}r^2)$$ satisfies the required condition, the last
term of (2.38) vanishes, and we are done.

Now, we address the second part of the theorem.

From (2.25), (2.28), we see that $\hat\nabla\lambda=0$ means
$$\delta\lambda=D\lambda,\hspace{5mm}D\lambda=\nabla
\lambda+\frac{i}{h}[r,\lambda].\leqno{(2.42)}$$ Since
$\lambda\in\mathcal W(\mathcal L)$, $\delta^{-1}\lambda=0$, and
(2.18), (2.42)  yield
$$\lambda=\lambda_{00}+\delta^{-1}D\lambda=a+\delta^{-1}D\lambda,
\leqno{(2.43)}$$ and, with the decomposition
$\lambda=\sum_{p=0}^\infty\hat\lambda_p$, we get the recurrence
formula
$$\hat\lambda_{p}=(\hat a)_p+\delta^{-1}\nabla
\hat\lambda_{p-1}+\frac{i}{h}\delta^{-1}([r,\lambda]_{p-1}),
\leqno{(2.44)}$$
which uniquely defines all the terms $\hat\lambda_p$.

Particularly, if $$a=\sum_{k=0}^\infty h^kf_k,\leqno{(2.45)}$$ the
first eight terms $\lambda_p$ are given by
$$\hat\lambda_p=\sum_{s=0}^{[p/2]}h^s(\delta^{-1}\nabla)^{p-2s}f_s,
\hspace{5mm}0\leq p\leq4,\leqno{(2.46)}$$
$$\hat\lambda_p=\sum_{s=0}^{[p/2]}h^s(\delta^{-1}\nabla)^{p-2s}f_s
\leqno{(2.47)}$$
$$+\frac{i}{h}\sum_{s=0}^{p-5}(\delta^{-1}\nabla)^s\delta^{-1}
\sum_{j=3}^{p-s-2}[\hat r_j,\hat\lambda_{p-s-j-1}]\hspace{5mm}5\leq p\leq8.$$

The terms $\hat\lambda_p$ for larger values of $p$ include
commutants with several factors $r$. If for any $\mu\in\hat
W(\mathcal L)$ we denote
$$\bar\mu=\sum_{s=0}^\infty(\delta^{-1}\nabla)^s\mu,
\leqno{(2.48)}$$ and reorder the terms of $\lambda$, as determined
by (2.46), (2.47), etc. we obtain
$$\lambda=\bar a+\frac{i}{h}\delta^{-1}[r,\bar a]
+(\frac{i}{h})^2\delta^{-1}[r,\delta^{-1} \overline{[r,\bar
a}]]+\cdots. \leqno{(2.49)}$$

Finally, we must check that (2.43) implies $\hat\nabla\lambda=0$.
First, we notice that (2.43) implies
$$\delta^{-1}\hat\nabla\lambda=\delta^{-1}D\lambda
-\delta^{-1}\delta\lambda=\lambda-a-\delta^{-1}
\delta\lambda=\delta\delta^{-1}\lambda=0.\leqno{(2.50)}$$
This allows us to use the Hodge decomposition (2.18) for $\hat\nabla\lambda$,
the Abelian character of $\hat\nabla$, and (2.43) to get
$$\hat\nabla\lambda=\delta^{-1}\delta\hat\nabla\lambda=
\delta^{-1}(D\hat\nabla\lambda-\hat\nabla^2\lambda)=
\delta^{-1}D\hat\nabla\lambda.\leqno{(2.51)}$$

Since the operator $\delta^{-1}D$ raises the degree, (2.51) yields
a recurrence relation for the homogeneous terms of
$\hat\nabla\lambda$, in the sense of the decomposition (2.3),
which shows that $\hat\nabla\lambda=0$. Q.e.d.

As a consequence of this main theorem we see that there exists
an injection
$$l:\mathcal C[[h]]
\rightarrow\mathcal W(\mathcal L),\leqno{(2.52)}$$ which sends the
formal power series $a\in\mathcal Z(\mathcal L))$ to the
$\hat\nabla$-parallel section $\lambda$ of $\mathcal W(\mathcal
L)$ which has the central projection $\lambda_{00}=a$.

This injection precisely is Fedosov quantization. The reason to
see it as a quantization process with links to quantum physics, is
that the mapping $l$ leads to a deformation of the commutative
product $fg$ into a non commutative product, also known as a {\it
star product} \cite{B}, namely,
$$f*g=l^{-1}(l(f)\circ l(g)),\hspace{5mm}f,g\in\mathcal C.\leqno{(2.53)}$$
The definition is correct since (2.21) and (2.28) show that
$l(f)\circ l(g)$ belongs to the image of $l$.

In various concrete cases such as the ones in Examples 1.5-1.9 the
mapping $l$ is defined by   corresponding versions of formulas
(2.37), (2.46), (2.47), (2.49), etc. In particular, since a Lie
algebroid always has connections, our presentation shows that
Fedosov quantization works for the regular Poisson manifolds and,
also, for the quasi-symplectic Poisson manifolds, even if the
latter are non regular. In other cases, the vanishing Atiyah class
condition is required, and it may be very restrictive. For
instance, this happens in the case of holomorphic symplectic
manifolds \cite{A}.
\vspace*{1cm}
{\small Department of Math.\\Univ. of Haifa\\Haifa 31905, Israel\\E-mail:
vaisman@math.haifa.ac.il}
\end{document}